\documentclass[12pt]{article}
\usepackage{a4, amsmath, amssymb, epsfig}

\newtheorem{thm}{Theorem}[section]
\newtheorem{lem}[thm]{Lemma}
\newtheorem{cor}[thm]{Corollary}
\newtheorem{prop}[thm]{Proposition}
\newtheorem{defi}[thm]{Definition}

\newtheorem{exam}[thm]{Example}

\newtheorem{question}[thm]{Question}

\newtheorem{observ}[thm]{Observation}

\newcounter{romcount}
\newenvironment{romlist}{\begin{list}{{\rm(\roman{romcount})}}{\usecounter%
{romcount}\setlength{\rightmargin}{\leftmargin}}}{\end{list}}
 
\newcommand{\ls}[2]{\: \makebox[0mm]{}^{#1}\! #2} 
\newcommand{\proof}{{\sc Proof.\ }} 
\newcommand{\qed}{\hspace*{\stretch{1}}$\blacksquare$ \\} 
\newcommand{\supp}{\operatorname{supp}} 
\newcommand{\nat}{\mathbb{N}}           
\newcommand{\Wre}[2]{\ensuremath{#1\,  ${\rm {Wr}}$\,  #2}} 
\newcommand{\wre}[2]{\ensuremath{#1\,  ${\rm {wr}}$\,  #2}} 

\def\a{\alpha}         
\def\b{\beta}          
\def\id{\equiv}            
\def\op{\sharp}            
\def\s{\sigma}         
\def\pad{\$}               
\def\dA{\delta _A}         

\newcommand{\ov}{\overline}

\begin{document}

\title{\large \bf Automatic structures for semigroup constructions}
\author{\normalsize I. Andrade, L. Descal\c{c}o and Manuel A. Martins \\
\normalsize Universidade de Aveiro,\\
\normalsize Departamento de Matem\'{a}tica,\\
\normalsize 3810-193 Aveiro, Portugal\\
\\
\normalsize Tel. +351234370359, Fax. +351234382014 \\
\normalsize e-mails: iandrade@oninet.pt, luisd@mat.ua.pt, martins@mat.ua.pt}
\maketitle

\begin{abstract}
We survey results concerning automatic structures for
semigroup constructions, providing references and
describing the corresponding automatic structures. The constructions we
consider are: free products, direct products, Rees matrix
semigroups, Bruck-Reilly extensions and wreath products.
\end{abstract}

\vspace{2em} \setlength{\arraycolsep}{0.3mm}


\section{Introduction}

The notion of ``automaticity'' has been widely studied in groups
(see \cite{AGA} and \cite{WPG} for example), and some progress has
been made in understanding the notion in the wider context of
semigroups. Many results about automatic semigroups concern
automaticity of standard semigroups constructions. We survey these
results for free products, direct products, Rees matrix
semigroups, Bruck-Reilly extensions and wreath products.

Some references on automatic semigroups are \cite{AutSem}
(introduction), \cite{SilvaGeo} (geometric aspects and
p-automaticity), \cite{otto2}, \cite{otto3}, \cite{otto1}
(computational and decidability aspects), \cite{DirProd}, \cite{RMAS},
\cite{AutMonCG} (semigroup constructions), \cite{Hoffmann}
(other notions of ``automaticity'' for semigroups) and
\cite{paperbimon}, \cite{ComNA} (examples).

We start by introducing the definitions we require. Given a non empty finite
set $A$, which we call an {\it alphabet}, we denote by $A^+$ the
free semigroup generated by $A$ consisting of finite sequences of
elements of $A$, which we call {\it words}, under the
concatenation; and  by $A^*$ the free monoid generated by $A$
consisting of $A^+$ together with the empty word $\epsilon$, the identity in $A^*$.
Let
$S$ be a semigroup and $\psi : A \rightarrow S$ a mapping. We say
that $A$ is a {\it finite generating set for $S$ with respect to
$\psi$} if the unique extension of $\psi$ to a semigroup
homomorphism $\psi: A^+ \rightarrow S$ is surjective. For $u,v \in
A^+$ we write $u \equiv v$ to mean that $u$ and $v$ are equal as
words and $u = v$ to mean that $u$ and $v$ represent the same
element in the semigroup i.e. that $u \psi = v \psi$. We say that
a subset $L$ of $A^*$, usually called a {\it language}, is {\it
regular} if there is a finite state automaton accepting $L$. To be
able to deal with automata that accept pairs of words and to
define automatic semigroups we need to define the set $A(2,\$)=
((A \cup\{\$\})\times(A\cup \{\$\})) \backslash \{(\$,\$)\}$ where
$\$$ is a symbol not in $A$ (called the {\em padding symbol}) and
the function $\delta_A: A^* \times A^* \rightarrow A(2,\$)^*$
defined by
\begin{equation*}
(a_1\ldots a_m,b_1\ldots b_n ) \delta_A = \left\{
\begin{array}{ll}
\epsilon & \mbox{ if } 0 = m = n \\
(a_1,b_1)\ldots(a_m,b_m) & \mbox{ if } 0 < m = n \\
(a_1,b_1)\ldots(a_m,b_m)(\$,b_{m+1})\ldots(\$,b_n) &
    \mbox{ if } 0 \le m < n \\
(a_1,b_1)\ldots(a_n,b_n)(a_{n+1},\$) \ldots (a_m,\$) &
    \mbox{ if } m > n \ge 0.
\end{array}
\right.
\end{equation*}

Let $S$ be a semigroup and $A$ a finite generating set for $S$
with respect to $\psi : A^+ \rightarrow S$. The pair $(A,L)$ is an
{\it automatic structure for $S$ (with respect to $\psi$)} if
\begin{itemize}
\item $L$ is a regular subset of $A^+$ and $L \psi = S$, \item
$L_= = \{(\alpha,\beta): \alpha,\beta \in L, \alpha = \beta\}
\delta_A \mbox{ is regular in } A(2,\$)^+, \mbox{ and}$ \item $L_a
= \{(\alpha,\beta):\alpha,\beta \in L, \alpha a = \beta \}
\delta_A \mbox{ is regular in } A(2,\$)^+ \mbox{ for each } a \in
A.$
\end{itemize}
We say that a semigroup is {\it automatic} if it has an automatic
structure.

We say that the pair $(A,L)$ is an {\it automatic structure with uniqueness} (with respect to $\psi$) for a semigroup $S$, if it
is an automatic structure and each element in $S$ is represented by an unique word in $L$ (the restriction of $\psi$ to $L$ is a bijection).
It is known (see \cite{AutSem}) that,
any automatic semigroup admits an
automatic structure with uniqueness.

We say that a semigroup is {\it prefix-automatic} or
{\it p-automatic} if it has an automatic
structure $(A,L)$ such that the set
$$
L_=^\prime = \{ (w_1,w_2) \delta_A: w_1 \in L, w_2 \in \mbox{{\rm Pref}}(L), w_1 = w_2 \}
$$
is also regular, where
$$\mbox{{\rm Pref}}(L) = \{w \in A^+: w w^\prime \in L \mbox{ for some }
w^\prime \in A^* \}.
$$

We will now present a result from \cite{FinReesInd} useful
to obtain automatic structures for the constructions considered in
the sections following.

We say that $T$ is a {\it subsemigroup of $S$ of finite Rees index} if
the set $S-T$ is finite.

\begin{prop}
\label{finite_rees_index} Let $S$ be a semigroup with a
subsemigroup $T$ of finite Rees index. Then $S$ is automatic if
and only if $T$ is automatic.
\end{prop}

From an automatic structure $(A,L)$ for $T$, an automatic
structure for $S$ can be easily obtained. Take $C$ to be a
finite set of new symbols in bijection with the elements of $S-T$,
$A^\prime = A \cup C$ and $L^\prime = L \cup C$. The pair
$(A^\prime,L^\prime)$ is an automatic structure for $S$.

The converse is not so trivial. We start from an automatic
structure with uniqueness $(A,L)$ for $S$. It was shown in
\cite{FinReesInd} that there exists a constant $k$ such that every
element of the set $\{ \a \in \mbox{sub}(L) : \ k \leq |\a| < 2k
\} $ maps to an element of $T$. The generating set $B$ for $S$
is
$$\{ b_\a : \a \in \mbox{sub}(L), k \leq |\a| < 2k \} \cup
\{ c_\a : \a \in L, |\a| < k, \a \in T \}$$
where the $b_{\a}$ and $c_{\a}$ are new symbols such that each
$b_\a$ and each $c_\a$ maps to the same element of $T$ as the
corresponding $\a$. The regular language $K$ is obtained as follows.
We take $U = \{ \a \in L : \a \mbox{ represents an
element of } T \}$ and let $\phi : U \to B^*$ be defined by
$$ \a \phi = \left\{ \begin{array}{ll}
        c_\a & \mbox{ for } |\a| < k \\
        b_{a_1 \dots a_k} b_{a_{k+1} \dots a_{2k} } \dots
b_{a_{(l-1)k+1} \dots a_{lk} }
  b_{a_{lk+1} \dots a_r} & \mbox{ for } |\a| \geq k
        \end{array} \right.
$$
with $k \leq r-lk < 2k$ and $\a \id a_1a_2\dots a_r$. Note that
 $\a \phi = \a$ in $T$. Taking $K = U
\phi$, the pair $(B, K)$ is an automatic structure for $T$.

In sections two, three and four we present results
about automaticy for free products, direct products and Rees
matrix semigroups, respectively. We omit the proofs and just
describe briefly how to obtain the corresponding automatic structures.
In sections five and six we present results from
\cite{LuisThesis} about automaticity of Bruck-Reilly extensions
and wreath products, respectively.

\section{Free Products}

If $S_1$ and $S_2$ are semigroups with presentations
$\langle A_1 | R_1 \rangle $ and $\langle A_2 | R_2 \rangle $ respectively
(where $A_1 \cap A_2 = \emptyset$),
then their free product $S_1 \ast S_2$ is the semigroup defined by
the presentation $\langle A_1 \cup A_2 | R_1 \cup R_2 \rangle $. The elements
of the product can bee seen as sequences $s_1\ldots s_m$ of elements of $S_1 \cup S_2$ such that two consecutive elements
do not belong to the same factor. The product of sequences $s_1\ldots s_m$, $s^\prime_1\ldots s^\prime_n$ is the concatenation
$s_1\ldots s_m s^\prime_1\ldots s^\prime_n$ if $s_m$ and $s^\prime_1$ do not belong to the same factor; otherwise it
is $s_1\ldots s_{m-1} s s^\prime_2\ldots s^\prime_n$, where $s$ is the product of $s_m$ by $s^\prime_1$ in their common factor.

Free products of semigroups and monoids were considered in \cite{AutSem}. For semigroups the following was
shown:

\begin{thm}
\label{free-products}
Let $S_1$ and $S_2$ be semigroups. Then $S_1 \ast S_2$ is automatic if
and only if both $S_1$ and $S_2$ are automatic.
\end{thm}

The proof of this theorem give us the automatic structures.
Suppose that $S_1$ and $S_2$ are automatic semigroups, with
automatic structures with uniqueness, say $(A_1, L_1)$ and $(A_2,
L_2)$ respectively, with $A_1$ and $A_2$ disjoint sets. Taking $A = A_1 \cup A_2$ and
$$L ~ = ~ (L_1 \cup \lbrace \epsilon \rbrace)(L_2 L_1)^\ast (L_2 \cup \lbrace \epsilon \rbrace) ~ - ~
\lbrace \epsilon \rbrace,$$
we obtain a pair $(A,L)$ which is an automatic structure for the
semigroup free product $S_1 * S_2$.
Conversely, suppose that $S_1 \ast S_2$ is automatic with an
automatic structure $(A, L)$. Letting
$$B = \lbrace a \in A : a \mbox{ represents an element of } S_1 \rbrace,$$
the pair $(B, L \cap B^+)$ is an automatic structure for $S_1$.

The monoid free product is not the same as the semigroup free product.
It is the same as the group free product and can be seen as the
semigroup free product with the identity subgroups amalgamated.
For monoids we have the following:

\begin{thm}
\label{free-products-2}
The monoid free product $M = M_1 * M_2$ is automatic if and only if both
monoids $M_1$ and $M_2$ are automatic.
\end{thm}

One implication was proved in \cite{AutSem}.
Suppose that $(A_1, L_1)$ and $(A_2, L_2)$ are automatic
structures with uniqueness for $M_1$ and $M_2$ respectively, with
$A_1 \cap A_2 = \lbrace e \rbrace$, $e$ representing the identity
element of each $M_i$, $e \in L_i$ $(i = 1, 2)$, and
$\overline{L}_i = L_i - \lbrace e \rbrace \subseteq (A_i - \lbrace
e \rbrace )^+$ $(i = 1, 2)$. Taking $A = A_1 \cup A_2$ and
$$L ~ = ~ \lbrace e \rbrace ( \overline{L}_1 \cup \lbrace \epsilon \rbrace )( \overline{L}_2 \overline{L}_1 )^\ast
(\overline{L}_2 \cup \lbrace \epsilon \rbrace ),$$
the pair $(A, L)$ is an automatic structure for $M$.

The converse was shown in \cite{AutMonCG}, answering a question
formalized in \cite{AutSem}. Finite generating sets $A_i$ for $M_i$, $i=1,2$, give us a finite
generating set $A= A_1 \cup A_2$ for $M$ with respect to an
homomorphism $\psi$. It is possible to obtain an automatic
structure $(A,L)$ for $M$ such that every element of $S$ is
represented by a unique element of $L$. The pair $(A,L(M_i))$,
where $L(M_i)=\{w\in L\,:\,w \psi \in M_i\}$, is an automatic
structure for $M_i$, $i=1,2$.

\section{Direct products}

The first result about direct products was obtained for monoids in
\cite{AutSem} where the authors have shown:

\begin{thm}
\label{Monoid-Direct-Product}
If $M_1$ and $M_2$ are automatic monoids, then their direct product
$M_1 \times M_2$ is automatic.
\end{thm}
Since we have the identities, an automatic structure for the
product can be obtained from automatic structures for the factors
in a natural way. We can start from automatic structures with
uniqueness $(A_1, L_1)$ and $(A_2, L_2)$ for $M_1$ and $M_2$
respectively, with $A_1 \cap A_2 = \emptyset$, $e_i \in A_i$,
$e_i$ representing the identity element of $M_i$, $e_i \in L_i$,
and $L_i - \lbrace e_i \rbrace \subseteq (A_i - \lbrace e_i
\rbrace )^+$ $(i = 1, 2)$. Let $\overline{L}_i$ denote $L_i -
\lbrace e_i \rbrace$ and let $A = A_1 \cup A_2$. For words $\a
\equiv a_1 \ldots a_n \in L_1$ and $\b \id b_1 \ldots b_m \in
L_2$, we define the word $\a \op \b \in A^+$ by
$$\a \op \b =
\left\{ \begin{array}{ll}
 a_1 b_1 \ldots a_n b_n
                                & \mbox{if } n = m, \\
 a_1 b_1 \ldots a_n b_n e_1 b_{n+1} \ldots e_1 b_m
                                & \mbox{if } n < m, \\
 a_1 b_1 \ldots a_m b_m a_{m+1} e_2 \ldots a_n e_2
                                & \mbox{if } n > m.
\end{array} \right.$$
If $\s : A(2, \pad)^\ast \to A^\ast$ is the homomorphism defined
by
$$(a, b) \mapsto a b, ~ ~ (a, \pad) \mapsto a e_2, ~ ~ (\pad, b)
\mapsto e_1 b,$$
then $\a \op \b = (\a, \b) \dA \s$. Let
$$L = \lbrace \a \op \b : \a \in L_1, \b \in L_2 \rbrace =
(L_1 \times L_2) \dA \s.$$
The pair $(A, L)$ is an automatic structure for $M = M_1 \times
M_2$.

Semigroups were then considered in \cite{DirProd} where the
authors have proved the following:
\begin{thm}
Let $S$ and $T$ be automatic semigroups.
\begin{romlist}
\item If $S$ and $T$ are infinite, then $S \times T$ is automatic
if and only if $S^2= S$ and $T^2 = T$. \item If $S$ is finite and
$T$ is infinite, then $S \times T$ is automatic if and only if
$S^2 = S$.
\end{romlist}
\label{autdirprod}
\end{thm}

In \cite{GRDP}, there were established necessary and
sufficient conditions for the direct product of semigroups to be
finitely generated:

\begin{prop}
Let $S$ and $T$ be two semigroups. If both $S$ and $T$ are
infinite then $S\times T$ is finitely generated if and only if
both $S$ and $T$ are finitely generated, $S^2 = S$ and $T^2 = T$.
If $S$ is finite and $T$ is infinite then $S \times T$ is finitely
generated if and only if $S^2 = S$ and $T$ is finitely generated.
\end{prop}

Using this result, Theorem \ref{autdirprod} has the
following equivalent formulation:

\begin{thm}
The direct product of automatic semigroups is automatic if and
only if it is finitely generated.
\end{thm}

The answer to the following converse question is not known even
for groups: If the direct product $G_1 \times G_2$ is automatic
are both factors $G_1$ and $G_2$ necessarily automatic?

Without the identities it is still possible to obtain automatic structures for
the product, starting from automatic structures for the factors, although the method is not
so natural.

For case (i) in Theorem \ref{autdirprod}, the general ideia is to start from two
automatic structures, say $(A,L)$ and $(B,K)$, for the factors $S$ and $T$ and then modify them,
using the fact that $S^2=S$ (and $T^2=T$) to control the length of the words in the languages,
in order to obtain new automatic structures. From the modified automatic structures,
say $(A^\prime,L^\prime)$ and $(B^\prime,K^\prime)$, an
automatic structure $(X,J)$ for the product $S \times T$ can be obtained by just taking
$X = A^\prime \times B^\prime$ and $J=\{(u_1,v_1),\ldots,(u_p,v_p): (u_i,v_i)\in X, u_1\ldots u_p \in L^\prime, v_1\ldots v_p \in K^\prime\}$.

For case (ii) we can assume that $S=\{s_1,\ldots,s_m\}$ and take an alphabet
$A=\{a_1,\ldots,a_m\}$ to represent the elements in $S$. Given an automatic structure
$(B,K)$ for $T$, the set $X = A \times B$ is a generating set for $S\times T$. Now, taking
$J=\{(u_1,v_1)\ldots(u_p,v_v): (u_i,v_i)\in X, v_1 \ldots v_p\in K\}$, the pair $(X,J)$ is an automatic
structure for the product $S \times T$ (the details can be found in \cite{DirProd}).

\section{Rees matrix semigroups}

The {\it Rees matrix semigroup} $S = \mathcal{M}[U;I,J;P]$ over the
semigroup $U$, with $P = (p_{j i})_{j \in J, i \in I}$ a matrix with
entries in $U$,
is the semigroup with the support set $I \times U
\times J$ and multiplication defined by $(l_1,s_1,r_1)(l_2,s_2,r_2)
= (l_1,s_1 p_{r_1 l_2} s_2, r_2)$ where $(l_1,s_1,r_1),(l_2,s_2,r_2)
\in I \times U \times J$.
We say that $U$ is the {\em base
semigroup} of the Rees matrix semigroup $S$.

We can obtain an automatic structure for a Rees matrix semigroup
$S$ by using the automatic structure for its base semigroup $U$,
as shown in \cite{RMAS}. We observe that the case where $U$ is a group, was
firstly considered in \cite{AutSimple}.
\begin{thm}
Let $S=\mathcal{M}[U;I,J;P]$ be a Rees matrix semigroup. If $U$ is
an automatic semigroup and if $S$ is finitely generated then $S$
is automatic.
\end{thm}

This theorem has the following equivalent formulation:

\begin{thm}
Let $S=\mathcal{M}[U;I,J;P]$ be a Rees matrix semigroup, where
$I,J$ are finite sets and $U \backslash V$ is finite, where $V$ is
the ideal of $U$ generated by the entries of the matrix $P$. If
$U$ is an automatic semigroup then $S$ is automatic. \label{thm1}
\end{thm}

In fact, it is described in \cite{RMAS} how to obtain an automatic
structure for the semigroup $S_1 = \mathcal{M}[U^1;I,J;P]$ from an
automatic structure with uniqueness for $V$ ($U^1$ stands for the monoid
obtained from $U$ by adding an identity). But note that, since
$S$ is finitely generated, $I$, $J$ and $U-V$ are finite, and so,
using Proposition \ref{finite_rees_index}, an automatic structure
for $S$ can then be obtained from an automatic structure for $U$.

We start from an automatic structure with uniqueness $(B,K)$ for
$V$, where $B = \{b_1,\ldots,b_n\}$ is a set of semigroup
generators for $V$. Then we write each $b_h\ (h \in N =
\{1,\ldots,n\})$ as $b_h = s_h p_{\rho_h \lambda_h} s_h^\prime$
where $s_h,s_h^\prime \in U^1, \rho_h \in J, \lambda_h \in I$. Let
$S_1 = \mathcal{M}[U^1;I,J;P]$. Given $(l,s,r) \in I \times V
\times J$ we can write $s = b_{\alpha_1}\ldots b_{\alpha_h}$ where
$b_{\alpha_1}\ldots b_{\alpha_h}$ is a word in $K$. So we can
write
$$
(l,s,r) = (l,s_{\alpha_1},\rho_{\alpha_1})
(\lambda_{\alpha_1},s_{\alpha_1}^\prime
s_{\alpha_2},\rho_{\alpha_2}) \ldots
(\lambda_{\alpha_h},s_{\alpha_h}^\prime,r).
$$
Since $U^1 \backslash V$ is finite and non empty we can write $U^1
\backslash V = \{x_1,\ldots,x_m\}$ with $m \ge 1$. We define a set
$A = C \cup D$ of semigroup generators for $S_1$ by
\begin{eqnarray}
& C = \{ c_{l i} : l \in I, i \in N \} \cup \{ d_{i j}: i,j \in N
\} \cup
\{ e_{j r} : j \in N, r \in J \},& \nonumber \\
 & D = \{f_{l h r}: l \in I, h \in \{1,\ldots,m\}, r \in J\} \nonumber &
\end{eqnarray}
with
$$
\begin{array}{lll}
\psi : A^+ \rightarrow S_1,\ & c_{l i}\mapsto (l,s_i,\rho_i), &
d_{i j} \mapsto (\lambda_i,s_i^\prime s_j, \rho_j),\\
&e_{j r} \mapsto  (\lambda_j,s_j^\prime,r), & f_{l h r} \mapsto
(l,x_h,r).
\end{array}
$$
Defining the language $L = L_1 \cup D$ to represent the elements
of $S_1$ with
$$
L_1 = \{ c_{l \alpha_1} d_{\alpha_1 \alpha_2} \ldots
d_{\alpha_{h-1} \alpha_h} e_{\alpha_h r}: b_{\alpha_1} \ldots
b_{\alpha_h} \in K, h \ge 1, l \in I, r \in J \},
$$
the pair $(A,L)$ is an automatic structure for $S_1$.

It was also shown in \cite{RMAS} that, in some particular
situations, it is possible to obtain an automatic structure for
the base semigroup, from the automatic structure for the
construction.

\begin{thm}
Let $S=\mathcal{M}[U;I,J;P]$ be a semigroup, and suppose that
there is an entry $p$ in the matrix $P$ such that $p U^1 = U$. If
$S$ is automatic then $U$ is automatic. \label{thm2}
\end{thm}

We start from an automatic structure with uniqueness $(A,L)$ for
the semigroup $S_1 = \mathcal{M}[U^1;I,J;P]$, where $A = \{ a_1,
\ldots, a_n \}$ is a generating set for $S_1$ with respect to
$$
\psi : A^+ \rightarrow S_1,\ a_h \mapsto (i_h,s_h,j_h)\ (h =
1,\ldots,n).
$$
The set $$B = \{ b_1,\ldots,b_n\} \cup \{c_{j i}: j \in J, i \in
I\}$$ is a generating set for $U^1$ with respect to
$$
\phi : B^+ \rightarrow U^1 ; b_h \mapsto s_h, c_{j i} \mapsto p_{j
i}\  (h = 1,\ldots,n, \ j \in J, i \in I).
$$
Without loss of generality we can assume that $p_{1 1} = p$. Let
$$L_{1 1} = L \cap (\{1\} \times U^1 \times \{1\}) \psi^{-1}.$$
Let
$$
f : A^+ \rightarrow B^+;
 a_{\alpha_1}  a_{\alpha_2} \ldots
a_{\alpha_h} \mapsto
 b_{\alpha_1}
c_{j_{\alpha_1} i_{\alpha_2}} b_{\alpha_2} \ldots
c_{j_{\alpha_{h-1}} i_{\alpha_h}} b_{\alpha_h} .
$$
Taking $K = L_{1 1} f$, the pair $(B,K)$ is an automatic structure
with uniqueness for $U^1$ with respect to $\phi$.

\begin{thm}
Let $S=\mathcal{M}[U;I,J;P]$ be a Rees matrix semigroup. If $S$ is
prefix-automatic then $U$ is automatic. \label{thm3}
\end{thm}

We start from a prefix-automatic structure with uniqueness $(A,L)$
for $S$ (see \cite{SilvaGeo}). We define $A,\psi,B,\phi,L_{1
1},f$ and $K$ as above just replacing $U^1$ by $U$ and
$S_1$ by $S$ in the definitions, and assume that $\psi
\restriction_A$ is injective. The pair $(B,K)$ is a
(prefix-)automatic structure with uniqueness for $U$ with respect
to $\phi$.

\section{Bruck-Reilly extensions}

Let $T$ be a monoid and $\theta: T \mapsto T$ be a monoid
homomorphism. The set $$\nat_0 \times T \times \nat_0$$ with the
operation defined by $$(m,t_1,n)(p,t_2,q) = (m-n+k, (t_1
\theta^{k-n})(t_2 \theta^{k-p}),q-p+k)\ (k = max\{n,p\}),$$ where
$\theta^0$ denotes the identity map on $M$, is called the
Bruck--Reilly extension of $T$ determined by $\theta$ and is
denoted by ${\rm BR}(T,\theta)$. The semigroup ${\rm
BR}(T,\theta)$ is a monoid with identity $(0,1_T,0)$, denoting by
$1_T$ the identity of $T$. This is a generalization of the
constructions from \cite{Bruck, Munn, Reilly}, also considered in
\cite{Isabel}.

\begin{thm}
If $T$ is a finite monoid, then any Bruck--Reilly extension of T
is automatic. \label{BRTfin}
\end{thm}
\proof Let $T = \{t_1,\ldots,t_l\}$ and let $\ov{T} =
\{\ov{t_1},\ldots,\ov{t_l}\}$ be an alphabet in bijection with
$T$. We define the alphabet $A = \{b,c\} \cup \ov{T}$ and the
regular language $$ L = \{c^m \ov{t} b^n: m,n \ge 0, \ov{t} \in
\ov{T}\} $$ on $A$. Defining the homomorphism $$ \psi: A^+
\rightarrow  {\rm BR}(T,\theta);\  \ov{t} \mapsto (0,t,0),\ c
\mapsto (1,1_T,0),\  b \mapsto (0,1_T,1) $$ it is clear that $A$
is a generating set for ${\rm BR}(T,\theta)$ with respect to
$\psi$ and, in fact, given an element $(m,t,n) \in \nat_0 \times T
\times \nat_0$, the unique word in $L$ representing it is $c^m
\ov{t} b^n$.

In order to prove that $(A,L)$ is an automatic structure with
uniqueness for $BR(T,\theta)$ we only have to prove that, for each
generator $a \in A$, the language $L_a$ is regular. To prove that
$L_b$ is regular we observe that $$(c^m \ov{t_i} b^n) b = (m, t_i,
n) (0,1_T,1) = (m,t_i,n+1) = c^m \ov{t_i} b^{n+1}$$ and so we can
write $$
\begin{array}{lll}
L_{b} & = & \displaystyle{\bigcup_{i=1}^{l}}  \{(c^m \ov{t_i}
b^n,c^m \ov{t_i} b^{n+1})\delta_A: n,m \in \nat_0\}\\ & = &
\displaystyle{\bigcup_{i=1}^{l}}(\{(c,c)\}^* \cdot
\{(\ov{t_i},\ov{t_i})\} \cdot \{(b,b)\}^* \cdot \{(\$,b)\})
\end{array}
$$ which is a finite union of regular languages and so is regular.
With respect to $L_c$ we have $$
\begin{array}{ll}
(c^m \ov{t}) c = (m,t,0)(1,1_T,0) = (m+1, t \theta,0) = c^{m+1}
\ov{t \theta},\\ (c^m \ov{t} b^{n+1}) c = (m, t, n+1) (1,1_T,0) =
(m, t, n) = c^m \ov{t} b^n\ (n,m \in \nat_0; \ov{t} \in \ov{T})
\end{array}
$$ and so we can write $$
\begin{array}{lll}
L_c & = & \displaystyle{\bigcup_{i = 1}^{l}} \{(c^m
\ov{t_i},c^{m+1} \ov{t_i \theta} )\delta_A: m \in \nat_0\} \cup \\
 & &  \displaystyle{\bigcup_{i = 1}^{l}} \{(c^m \ov{t_i} b^{n+1}, c^m
 \ov{t_i} b^n)\delta_A : m,n \in \nat_0 \}\\
& = & \displaystyle{\bigcup_{i=1}^{l}} (\{(c,c)\}^* \cdot
\{(\ov{t_i},c)(\$, \ov{t_i \theta})\}) \cup \\ & &
\displaystyle{\bigcup_{i=1}^{l}} (\{c,c)\}^* \cdot
\{(\ov{t_i},\ov{t_i})\} \cdot \{(b,b)\}^* \cdot \{(b,\$)\})
\end{array}
$$ and we conclude that $L_c$ is a regular language as well.

We now fix an arbitrary $\ov{t} \in \ov{T}$ and prove that
$L_{\ov{t}}$ is regular. For any words $c^m \ov{t_{\alpha}} b^n,
c^p \ov{t_{\beta}}b^q \in L$ we have $$ c^m \ov{t_{\alpha}} b^n
\ov{t} = c^p \ov{t_{\beta}}b^q $$ if and only if $m = p, n = q$,
and $t_{\alpha} (t \theta^n) = t_{\beta}$, because
$$ c^m
\ov{t_{\alpha}} b^n \ov{t} = (m,t_{\alpha},n)(0,t,0) = (m,
t_{\alpha} (t \theta^n), n). $$ Since $T$ is finite the set $\{ t
\theta^n: n \in \nat_0\}$ is finite as well. Taking $j$ to be
minimum such that the set $C_j = \{k \ge j: t \theta^j = t
\theta^{k+1}\}$ is non empty and $k$ to be the minimum element of
$C_j$, we will now show that
$$\{ t \theta^n: n \in \nat_0\} = \{t,
t \theta,\ldots, t \theta^j,\ldots, t \theta^{k}\}.$$ Given $n \ge
j$ we have $n = j + h$ with $h \ge 0$ and, dividing $h$ by
$k+1-j$, we obtain $n = j + q(k+1-j)+r$ with $q \ge 0$ and $0 \le
r < k+1 -j$. We now prove, by induction on $q$, that $t \theta^{j
+r + q(k+1-j)} = t \theta^{j+r}$ for $q \ge 0$. For $q=0$ it holds
trivially and for $q > 0$ we have
$$
\begin{array}{ll}
t \theta^{j +r + q(k+1-j)} & = t \theta^{j+r+k+1-j+(q-1)(k+1-j)}
= (t \theta^r)(t \theta^{k+1})(t \theta^{(q-1)(k+1-j)}) \\
 & = (t \theta^r)(t \theta^j)(t \theta^{(q-1)(k+1-j)}) =
t \theta^{j+r+(q-1)(k+1-j)}.
\end{array}
$$
We can then write
$$
\begin{array}{lll}
L_{\ov{t}} & = &  \displaystyle{\bigcup_{n=0}^{j-1}} \{(c^m
\ov{t_{\alpha}} b^n, c^m \ov{t_{\alpha} (t \theta^n)} b^n)\delta_A
:m \in \nat_0, t_{\alpha} \in T\} \cup \\
& & \displaystyle{\bigcup_{n=j}^k} \{(c^m \ov{t_{\alpha}}
b^{n+q(k+1-j)}, c^m \ov{t_{\alpha} (t \theta^n)}
b^{n+q(k+1-j)})\delta_A: m,q \in \nat_0, t_{\alpha} \in T \} \\
& = & \displaystyle{\bigcup_{n=0}^{j-1}} (\{(c,c)\}^* \cdot
\{(\ov{t_{\alpha}},\ov{t_{\alpha} (t \theta^n)}): t_{\alpha} \in
T\} \cdot \{(b,b)\}^*) \cup\\
&  & \displaystyle{\bigcup_{n=j}^k} (\{(c,c)\}^* \cdot
\{(\ov{t_{\alpha}},\ov{t_{\alpha} (t \theta^n)}): t_{\alpha} \in
T\} \cdot \{(b,b)^n\} \cdot \{(b,b)^{k+1-j}\}^*)
\end{array}
$$ and since all sets in this union are regular we conclude that
$L_{\ov{t}}$ is regular as well. \qed

From now on we assume that $T$ is an automatic monoid and we fix
an automatic structure $(X,K)$ with uniqueness for $T$, where $X =
\{x_1,\ldots,x_n\}$ is a set of semigroup generators for $T$ with
respect to the homomorphism
$$\phi : X^+ \rightarrow T.$$
We define the alphabet
\begin{equation}
A = \{b,c\} \cup X \label{eqA}
\end{equation}
to be a set of semigroup generators for ${\rm BR}(T,\theta)$ with
respect to the homomorphism $$ \psi: A^+  \rightarrow  {\rm
BR}(T,\theta),  x_i \mapsto (0,x_i \phi,0), c \mapsto (1,1_T,0), b
\mapsto (0,1_T,1), $$ and the regular language
\begin{equation}
L = \{ c^i w b^j: w \in K; i,j \in \nat_0\} \label{eqL}
\end{equation}
on $A^+$, which is a set of unique normal forms for
 ${\rm BR}(T,\theta)$, since we have $(c^i w b^j)\psi = (i, w\phi,j)$
for $w \in K$, $i,j \in \nat_0$. As usual, to simplify notation,
we will avoid explicit use of the homomorphisms $\psi$ and $\phi$,
associated with the generating sets, and it will be clear from the
context whenever a word $w \in X^+$ is being identified with an
element of $T$, with an element of ${\rm BR}(T,\theta)$ or
considered as a word. In particular, for a word $w \in X^+$ we
write $w \theta$ instead of $(w \phi) \theta$, seeing $\theta$
also as a homomorphism $\theta: X^+ \rightarrow T$, and we will
often write $(i,w,j)$ instead of $(i,w \phi,j)$ for $i,j \in
\nat_0$.

For $(A,L)$ to be an automatic structure for ${\rm BR}(T,\theta)$
the languages $$
\begin{array}{lll}
L_b = & \{(c^i w b^j,c^i w b^{j+1})\delta_A: w \in K; i,j \in
\nat_0\}, \\ L_c = & \{(c^i w b^{j+1},c^i w b^j)\delta_A: w \in K;
i,j \in \nat_0\} \cup \\ & \{(c^i w_1,c^{i+1} w_2)\delta_A:
w_1,w_2 \in K; i \in \nat_0; w_2 = w_1 \theta \}, \\ L_{x_r} = &
\{(c^i w_1 b^j,c^i w_2 b^j) \delta_A: (w_1,w_2)\delta_X \in K_{x_r
\theta^j}; i,j \in \nat_0 \}\ (x_r \in X),
\end{array}
$$ must be regular. The language $L_b$ is regular, since we have
$$L_b = \{(c,c)\}^* \cdot \{(w,w)\delta_X: w \in K\} \cdot
\{(b,b)\}^* \cdot \{(\$,b)\},$$ but there is no obvious reason why
the languages $L_c$ and $L_{x_r}$ should also be regular. We will
consider particular situations where $(A,L)$ is an automatic
structure for ${\rm BR}(T,\theta)$.

\begin{thm}
If $T$ is an automatic monoid and $\theta: T \rightarrow T;\ t
\mapsto 1_T$ then ${\rm BR}(T,\theta)$ is automatic.
\label{thmmapsto1}
\end{thm}

To show this we use the notion of {\it padded product of languages} and an
auxiliary result whose proof can be found in \cite{LuisThesis}.
Fixing an alphabet $A$, and given two regular languages $M,N$ in $(A^* \times A^*)\delta$,
the {\it padded product of languages} $M$ and $N$ is
$$
M \odot N = \{ (w_1 w_1^\prime, w_2 w_2^\prime)\delta:
(w_1,w_2) \delta \in M, (w_1^\prime, w_2^\prime)\delta \in N\}
$$
The result is the following:

\begin{lem}
Let $A$ be an alphabet and let $M,N$ be regular languages on
$(A^* \times A^*)\delta$.
If there exists a constant $C$ such that, for any two
words $w_1, w_2 \in A^*$ we have
$$
(w_1,w_2) \delta \in M \; \Rightarrow \; | |w_1| - |w_2|| \le C,
$$
then the language $M \odot N$ is regular.
\label{thmconcat}
\end{lem}

{{\sc Proof. of Theorem \ref{thmmapsto1}\ }} To show that the pair $(A,L)$ defined by (\ref{eqA}) and
(\ref{eqL}) is an automatic structure for ${\rm BR}(T,\theta)$ we
just have to prove that the languages $L_c$ and $L_{x_r}\ (x_r \in
X)$ are regular. But now, denoting by $w_{1_T}$ the unique word in
$K$ representing $1_T$, we have
$$
\begin{array}{lll}
L_c & = &  \{(c^i w b^{j+1},c^i w b^j)\delta_A: w \in K; i,j \in
\nat_0\}
 \cup \{(c^i w,c^{i+1} w_{1_T})\delta_A: w \in K; i \in \nat_0 \} \\
& = & (\{(c,c)\}^* \cdot \{(w,w)\delta_X: w \in K\} \cdot
\{(b,b)\}^* \cdot \{(b,\$)\}) \cup \\ & & ((\{(c,c)\}^* \cdot
\{(\$,c)\}) \odot (K \times \{w_{1_t}\})\delta_X ),
\end{array}
$$
which is a regular language by Lemma
\ref{thmconcat}. We have
$$
\begin{array}{lll}
L_{x_r} & = & \{(c^i w b^j,c^i w b^j) \delta_A: w \in K,i \in
\nat_0, j\in \nat \} \cup \\ & & \{(c^i w_1,c^i w_2)\delta_A :
(w_1,w_2)\delta_X \in K_{x_r}; i \in \nat_0 \} \\ & = &
(\{(c,c)\}^* \cdot \{(w,w)\delta_X: w \in K\} \cdot \{(b,b)\}^+)
\cup \\ & & (\{(c,c)\}^* \cdot K_{x_r})
\end{array}
$$ because, for any $c^i w b^j \in L$ with $j \ge 1$, we have
$$(c^i w b^j) x_r = (i,w,j)(0,x_r,0) = (i,w (x_r \theta^j), j) =
(i, w, j) = c^i w b^j$$ and for $c^i w \in L$  we have $$(c^i w)
x_r = (i,w,0)(0,x_r,0) = (i,w x_r,0).$$ Therefore $L_{x_r}$ is
also a regular language and so ${\rm BR}(T,\theta)$ is automatic.
\qed

\begin{thm}
If $T$ is an automatic monoid and $\theta$ is the identity in $T$
then ${\rm BR}(T,\theta)$ is automatic.
\end{thm}
\proof We use the generating set $A$ defined by equation
(\ref{eqA}) but we now define $L = \{c^i b^j w : w \in K \}$
observing that, since $\theta$ is the identity, for any $x_r \in
X$, we have $$
\begin{array}{ll}
x_r c = (0,x_r,0)(1,1_T,0) = (1, x_r \theta, 0)=
   (1, x_r, 0) = (1, 1_T, 0)(0, x_r, 0) = c x_r,\\
x_r b = (0, x_r, 0) (0, 1_T, 1) = (0, x_r, 1) =
   (0, x_r \theta, 1) = (0, 1_T, 1) (0, x_r, 0) = b x_r.
\end{array}
$$
The language $L$ is regular and it is a set of unique normal forms
for ${\rm BR}(T,\theta)$. Also the languages $$
\begin{array}{lll}
L_b & = & \{(c^i b^j w, c^i b^{j+1} w)\delta_A: w \in K; i,j \in
\nat_0\} \\
 & = & (\{(c,c)\}^* \cdot \{(b,b)\}^* \cdot \{(\$,b)\}) \odot
       \{(w,w)\delta_X: w \in K\}, \\
L_c & = & \{(c^i b^{j+1} w, c^i b^j w)\delta_A: w \in K; i,j \in
\nat_0\} \cup \\ & & \{(c^i w,c^{i+1} w)\delta_A: i \in \nat_0,
w\in K\} \\ & = & ((\{(c,c)\}^* \cdot \{(b,b)\}^* \cdot
\{(b,\$)\}) \odot \{(w,w)\delta_X: w \in K\}) \cup \\ & &
((\{(c,c)\}^* \cdot  \{(\$,c)\}) \odot \{(w,w)\delta_X: w \in
K\}),\\ L_{x_r} & = & \{(c^i b^j w_1,c^i b^j w_2)\delta_A:
 (w_1,w_2) \delta_X \in K_{x_r} \} \\
 & = & (\{(c,c)\}^* \cdot \{(b,b)\}^*) \cdot K_{x_r}
\end{array}
$$ are regular, by Lemma \ref{thmconcat}, and so $(A,L)$ is an
automatic structure for ${\rm BR}(T,\theta)$. \qed

We say that a semigroup $T$ is of {\it finite geometrical type}
(fgt) (see \cite{SilvaGeo}) if for every $t_1 \in T$, there
exists $k \in \nat$ such that the equation
$$x t_1 = t_2$$
has at most $k$ solutions for every $t_2 \in M$.

To prove next theorem we will use the following two auxiliary results from \cite{LuisThesis}:

\begin{lem}
Let $T$ be a {\rm fgt} monoid with an automatic structure with
uniqueness $(X,K)$. Then for every $w \in X^+$ there is a constant
$C$ such that $(w_1,w_2)\delta_X \in K_w$ implies $||w_1| -|w_2||
< C$. \label{lemfgt}
\end{lem}

\begin{lem}
\label{FiniteReg}
Let $S$ be a finite semigroup, $X$ be a finite set and $\psi: X^+
\rightarrow S$ be a surjective homomorphism. For any $s \in S$ the
set $s \psi^{-1}$ is a regular language.
\end{lem}

\begin{thm}
Let $T$ be a {\it fgt} automatic monoid and let $\theta: T
\rightarrow T$ be a monoid homomorphism. If $T \theta$ is finite
then ${\rm BR}(T,\theta)$ is automatic.
\end{thm}
\proof We will prove that the pair $(A,L)$ defined by (\ref{eqA})
and (\ref{eqL}) is an automatic structure for ${\rm
BR}(T,\theta)$. We have $$
\begin{array}{ll}
L_c =  & \{(c^i w b^{j+1},c^i w b^j)\delta_A: w \in K; i,j \in
\nat_0\} \cup \\
 & \{(c^i w_1,c^{i+1} w_2)\delta_A: w_1,w_2 \in K;
i \in \nat_0; w_2 = w_1 \theta \}
\end{array}
$$
and, since the language
$$
\begin{array}{l}
\{(c^i w b^{j+1},c^i w b^j)\delta_A: w
\in K; i,j \in \nat_0\} =\\
 \{(c,c)^i\}^* \cdot \{(w,w)\delta_X: w
\in K\} \cdot \{(b,b)\}^* \cdot \{(b,\$)\}
\end{array}
$$
is regular, we just have to prove that the language
$$M = \{(c^i w_1,c^{i+1}
w_2)\delta_A: w_1,w_2 \in K; i \in \nat_0; w_2 = w_1 \theta \}$$
is also regular. For any $t \in T \theta$ let $w_t$ be the unique
word in $K$ representing $t$. Let $$
\begin{array}{lll}
N = & \{(w_1,w_2)\delta_X: w_1,w_2 \in K; w_2 = w_1 \theta \} = \\
& \displaystyle{\bigcup_{t \in T \theta}} \{(w_1,w_2)\delta_X:
w_1,w_2 \in K; w_2 = w_1 \theta = t\} = \\ &
\displaystyle{\bigcup_{t \in T \theta}} \{(w_1,w_t)\delta_X: w_1
\in K; w_1 \in
 (t \theta^{-1}) \phi^{-1} \} = \\
& \displaystyle{\bigcup_{t \in T \theta}} ( ((t \theta^{-1})
\phi^{-1} \cap K ) \times \{w_t\}) \delta_X.
\end{array}
$$
We can define $\psi: X^+ \rightarrow T \theta;\ w \mapsto w \phi
\theta$ and, since $T \theta$ is finite, for any $t \in T \theta$,
we can apply Lemma \ref{FiniteReg} and conclude that $(t
\theta^{-1}) \phi^{-1} = t \psi^{-1}$ is regular. Therefore, $N$ is
a regular language and, since we have $$
\begin{array}{ll}
M = & \{(c^i w_1,c^{i+1} w_2)\delta_A: (w_1,w_2)\delta_X \in N; i
\in \nat_0\} = \\ & (\{(c,c)\}^* \cdot \{(\$,c)\}) \odot N,
\end{array}
$$ by Lemma \ref{thmconcat}, $M$ is a regular language as well. We
will now prove that the language $$L_{x_r} = \{(c^i w_1 b^j,c^i
w_2 b^j) \delta_A: (w_1,w_2)\delta_X \in K_{x_r \theta^j}; i,j \in
\nat_0 \}$$ is regular. Since $T \theta$ is finite we can, as in
the proof of Lemma \ref{BRTfin}, take $j,k$ to be minimum with
$x_r \theta^j = x_r \theta^{k+1}$ and $j \le k$, and we have $x_r
\theta^{j+r+q(k+1-j)} = x_r \theta^{j+r}$ for $j \le j+r < k+1$
and $q \ge 0$. Therefore, we can write
$$
\begin{array}{llll}
L_{x_r} & = & \displaystyle{\bigcup_{n=0}^{j-1}} & \{(c^i w_1
b^n,c^i w_2 b^n)\delta_A:
              (w_1,w_2)\delta_X \in K_{x_r \theta^n}; i \in \nat_0\} \cup\\
        & & \displaystyle{\bigcup_{n=j}^k} & \{(c^i w_1 b^{n+q(k+1-j)},
          c^i w_2 b^{n+q(k+1-j)})\delta_A : (w_1,w_2)\delta_X \in K_{x_r
          \theta^n}; i,q \in \nat_0\}\\
& = & \displaystyle{\bigcup_{n=0}^{j-1}} & (\{(c,c)\}^* \cdot
(K_{x_r \theta^n} \odot \{(b,b)\}^*)) \cup\\ &  &
\displaystyle{\bigcup_{n=j}^k}& (\{(c,c)\}^* \cdot (K_{x_r
\theta^n} \odot (\{(b,b)^n\} \cdot \{(b,b)^{k+1-j}\}^*))).
\end{array}
$$
Since $T$ is fgt, by Lemma \ref{lemfgt} there is a constant $C$
such that
$$(w_1,w_2)\delta_X \in K_{x_r \theta^n} \; \Rightarrow \;
||w_1|-|w_2|| < C$$ for any $n =0,\ldots,k$, and therefore we can
apply Lemma \ref{thmconcat} and we conclude that $L_{x_r}$ is a
regular language. \qed

Since automatic groups are characterized by the fellow traveller
property and Bruck--Reilly extensions of groups are somehow
``almost groups'' the following is a natural question:
Is a Bruck--Reilly extension of a group automatic if and only if
it has the fellow traveller property?

\section{Wreath products}

We consider the automaticity of the wreath product of semigroups,
$\wre{S}{T}$, in the case where $T$ is a finite semigroup. We
start by giving the necessary and sufficient conditions, obtained
in \cite{FGWP}, for the wreath product, to be finitely
generated, when $T$ is finite. Finite generation of the wreath product is related to
finite generation of the diagonal $S$-act. We use the conditions
obtained for the case where the diagonal $S$-act is not finitely
generated to prove that, in this case, the wreath product
$\wre{S}{T}$ is automatic whenever it is finitely generated and
$S$ is an automatic semigroup.

We start by giving the definitions we require. If $S$ is a
semigroup and $X$ is a set, then the set $S^X$ of all mappings $X
\rightarrow S$ forms a semigroup under {\it component-wise}
multiplication of mappings: for $f,g \in S^X$, $fg: X \rightarrow
S;\ x \mapsto (xf)(xg)$; this semigroup is called the {\it
Cartesian power} of $S$ by $X$. If $S$ has a distinguished
idempotent $e$, then the {\it support} of $f \in S^X$ relative to
$e$ is defined by $$ {\rm supp}_e(f)= \{x \in X: x f \neq e\}. $$
The set $$ S^{(X)_e} = \{f \in S^X: |\supp_e( f)| < \infty \} $$
is a subsemigroup of $S^X$; it is called the {\it direct power} of
$S$ relative to $e$ (the subscript $e$ is usually omitted).
If $X$ is finite of size $n$
then $S^X$ and $S^{(X)_e}$ coincide, and they are isomorphic to
the semigroup $S^{(n)}$ consisting of $n$-tuples of elements of
$S$ under the component-wise multiplication. In this context, we
write $S^{(X)_e}$ even if $S$ has no idempotents; we can think of
this as computing supports with respect to an identity adjoined to
$S$.

The {\it unrestricted wreath product} $\Wre{S }{ T}$ of two
semigroups is the set $S^T \times T$ under multiplication $$
(f,t)(g,u) = (f \ls{t}{g}, t u), $$ where $\ls{t}{g} \in S^T$ is
defined by $$ (x) \ls{t}{g} = (x t) g. $$ Let $e \in S$ be a
distinguished idempotent. The {\it (restricted) wreath product}
$\wre{S_e}{T}$ (with respect to $e$) is the subsemigroup of
$\Wre{S}{T}$ generated by the set $\{(f,t) \in \Wre{S}{T}: |
\supp_e(f)| < \infty \}$ (again the subscript $e$ is often
omitted).

The wreath product $\wre{S}{T}$ coincides with the unrestricted
wreath product $\Wre{S}{T}$ in the case where $T$ is finite, as
observed in \cite[Chapter 3]{Rob}.

An {\it action} of a semigroup $S$ on a set $X$ is a mapping $X
\times S \rightarrow X$, $(x,s) \mapsto x s$, satisfying $(x
s_1)s_2 = x (s_1 s_2)$. The set $X$, together with an action, is
called an {\it S-act}. It is said to be generated by a set $U
\subseteq X$ if $U S^1 = X$, and {\it finitely generated} if there
exists a finite such $U$.

The {\it diagonal act} of a semigroup $S$ is the set $S\times S$
with the action $(s_1,s_2)s = (s_1 s,s_2 s)$. The
diagonal acts of infinite groups, free semigroups, free
commutative semigroups and completely simple semigroups are not
finitely generated. On the other hand, the diagonal act of the full transformation
monoid $T_\nat$ on positive integers can be generated by a single
element; see \cite{Bulman}. In \cite{DAMon} the authors give an
example of an infinite, finitely presented monoid with a finitely
generated diagonal act.

We will only state the conditions obtained in \cite{FGWP} for the
case where $T$ is finite and $S$ is infinite.

\begin{prop}
Let $S$ be an infinite semigroup and let $T$ be a finite
non-trivial semigroup. If the diagonal $S$-act is finitely
generated then $\wre{S}{T}$ is finitely generated if and only if
the following conditions are satisfied:
\begin{romlist}
\item $S^2 = S$ and $T^2 = T$; \item $S$ is finitely generated.
\end{romlist}
If the diagonal $S$-act is not finitely generated then
$\wre{S}{T}$ is finitely generated if and only if the following
conditions are satisfied:
\begin{romlist}
\item $S^2 = S$; \item $S$ is finitely generated; \item every
element of $T$ is contained in the principal right ideal generated
by a right identity.
\end{romlist}
\label{proprob}
\end{prop}

We will now consider the automaticity of the wreath product
$\wre{S}{T}$ in the case where $T$ is finite. In the case where
$S$ is also finite, $\wre{S}{T}$ is finite as well, and, in
particular, it is automatic. We will consider the case where $S$
is infinite and the diagonal $S$-act is not finitely generated.

\begin{thm}
If $S$ and $T$ are semigroups satisfying the following conditions:
\begin{romlist}
\item $T$ is finite; \item $S$ is automatic; \item the diagonal
$S$-act is not finitely generated; \item the wreath product
$\wre{S}{T}$ is finitely generated;
\end{romlist}
then $\wre{S}{T}$ is automatic. \label{thmwr}
\end{thm}

To prove this theorem we will need some notation and a result from
\cite{RMAS}. A {\it generalized sequential machine} (gsm for
short) is a six-tuple $\mathcal{A} = (Q,A,B,$ $\mu,q_0,T)$ where
$Q$, $A$ and $B$ are finite sets, (called the {\it states}, the
{\it input alphabet} and the {\it output alphabet} respectively),
$\mu$ is a (partial) function from $Q \times A$ to finite subsets
of $Q \times B^+$, $q_0\in Q$ is the {\it initial state} and
$T\subseteq Q$ is the {\it set of terminal states}. We can read
$(q^\prime,u)\in (q,a)\mu $ in the following way: if $\mathcal{A}$ is in state $q$ and receives input $a$,
then it can move into state $q^\prime$ and output $u$.

We can interpret $\mathcal{A}$ as a directed labelled graph with vertices $Q$,
and an edge $q\xrightarrow{(a,u)} q^\prime$ for every
pair $ (q^\prime,u)\in (q,a)\mu$.
For a path
$$
\pi\::\: q_1\xrightarrow{(a_1,u_1)} q_2\xrightarrow{(a_2,u_2)} q_3\ldots
\xrightarrow{(a_n,u_n)} q_{n+1}
$$
we define
$$
\Phi(\pi)=a_1a_2\ldots a_n,\ \Sigma(\pi)=u_1u_2\ldots u_n.
$$
For $q,q^\prime\in Q$, $u\in A^+$ and $v\in B^+$ we write
$q\xrightarrow{(u,v)}_+ q^\prime$ to mean that there exists
a path $\pi$ from $q$ to $q^\prime$ such that $\Phi(\pi)\equiv u$
and $\Sigma(\pi)\equiv v$, and we say that $(u,v)$ is the
{\it label} of the path. We say that a path is {\it successful}
if it has the form $q\xrightarrow{(u,v)}_+ t$ with $t \in T$.

The gsm $\mathcal{A}$ induces a mapping
$\eta_{\mathcal{A}}\::\: {\mathcal{P}}(A^+)\longrightarrow
{\mathcal{P}}(B^+)$ from subsets of $A^+$ into subsets of $B^+$
defined by
$$
X\eta_{\mathcal{A}}=\{ v\in B^+\::\: (\exists u\in X)(\exists t\in T)
(q_0\xrightarrow{(u,v)}_+ t)\}.
$$
It is well known that if $X$ is regular then so is $X\eta_{\mathcal{A}}$;
see \cite{ForLang}.
Similarly, $\mathcal{A}$ induces a mapping $\zeta_{\mathcal{A}}\::\:
{\mathcal{P}}(A^+\times A^+)\longrightarrow {\mathcal{P}}(B^+\times B^+)$
defined by
$$
Y\zeta_{\mathcal{A}}=\{ (w,z)\in B^+\times B^+\::\:
(\exists (u,v)\in Y)(w\in u\eta_{\mathcal{A}}\ \&\ z\in v\eta_{\mathcal{A}})\}.
$$

The next lemma asserts that, under certain conditions, this
mapping also preserves regularity.

\begin{lem}
Let
$\mathcal{A} = (Q,A,B,\mu,q_0,T)$ be a {\rm gsm},
and let $\pi_A\::\: (A^\ast\times A^\ast)\delta_A\longrightarrow
A^\ast\times A^\ast$ be the inverse of $\delta_A$.
Suppose that there is a constant $C$ such that for any two paths
$\alpha_1, \alpha_2$ in $\mathcal{A}$, we have
\begin{equation}
| \Phi(\alpha_1) | = | \Phi(\alpha_2)| \Longrightarrow
||\Sigma(\alpha_1)| - | \Sigma(\alpha_2)|| \le C.\label{pre}
\end{equation}
If $M\subseteq (A^+\times A^+)\delta_A$ is a regular language in $A(2,\$)^+$
then
$N = M \pi_A \zeta_{\mathcal{A}} \delta_B$ is a regular
language in $B(2,\$)^+$. \label{thmgsm}
\end{lem}

Also the following simple fact, from \cite{LuisThesis}, will be
used in our proof.

\begin{lem}
\label{thmsses}
Let $S$ be an automatic semigroup such that $S^2 = S$. Then $S$
has an automatic structure with uniqueness
$(A,K)$ such that $K \cap A = \emptyset$.
\end{lem}

{{\sc Proof of Theorem \ref{thmwr}.}}  We assume, without loss of generality,
that $T = \{t_1,\ldots,t_m\}$ with $m > 1$. By
using Proposition \ref{proprob} we know that $S$ is finitely
generated and $S^2 = S$. So, by Theorem \ref{autdirprod}, we
conclude that the direct product $S^{|T|}$ is automatic. Let
$(F,K)$ be an automatic structure for $S^{|T|}$ with uniqueness
with $F = \{f_1,\ldots,f_k\}$. Since $S^2 =S$, we can use Lemma
\ref{thmsses}, and assume that $K$ does not have words of length
1. Given $t \in T$, using again Proposition \ref{proprob}, there
is a right identity $e \in T$ such that $t = e q$ for some $q \in
T$. So we can define a generating set
$$Y = \{e_1,\ldots,e_m\} \cup \{q_1,\ldots,q_m\}$$
for $T$ such that $t_i = e_i q_i$ for $i=1,\ldots,m$ and
$e_1,\ldots,e_m$ represent (not necessarily distinct) right
identities in $T$. We define a new alphabet $A$ by $$A =
\{(f,e_i): f \in F, i = 1, \ldots, m\} \cup \{(f,q_i):f \in F,
i=1,\ldots,m\}$$ and a language $L$ on $A$ by $$L =
\bigcup_{i=1,\ldots,m}{\{(f_{\alpha_1},e_i)\ldots(f_{\alpha_{n-1}},e_i)
(f_{\alpha_n},q_i): f_{\alpha_1} \ldots f_{\alpha_n} \in K\}}.$$
We will prove that the pair $(A,L)$ is an automatic structure for
$\wre{S}{T}$ (with uniqueness). To see that $A$ generates
$\wre{S}{T}$ and that $L$ is a set of unique representatives for
$\wre{S}{T}$ we observe that, given $(f,t_i) \in \wre{S}{T}$ there
is only one word $f_{\alpha_1}\ldots f_{\alpha_n}$ in $K$ such
that $f = f_{\alpha_1}\ldots f_{\alpha_n}$. So there is only one
word in $L$ representing $(f,t_i)$ which is
$$(f_{\alpha_1},e_i)\ldots(f_{\alpha_{n-1}},e_i)
(f_{\alpha_n},q_i).$$ To prove that $L$ is a regular language we
now define a gsm $\mathcal{A}$ such that $K \eta_{\mathcal{A}} =
L$. Let
$$\mathcal{A} = (Q,F,A,\mu,q_0,\{\chi\})$$
with $Q = \{q_0,\ldots,q_m\} \cup \{\chi\}$, where $q_0$ is the
initial state, $\chi$ is the only final state and $\mu$ is a
partial function from $Q \times F$ to finite subsets of $Q \times
A^+$ defined by: $$
\begin{array}{rcl}
(q_0,f) \mu & = & \{(q_i,(f,e_i)) \}\ (i = 1,\ldots,m), \\ (q_i,f)
\mu  & = & \{(q_i,(f,e_i)),(\chi,(f,q_i))) \}\ (i=1,\ldots,m ).
\end{array}
$$
We will now prove that $L_{(f,e_r)}$ is a regular language, for
$(f,e_r)\in A$. If we define
$$L_{(f,e_r)}^{(i)} = L_{(f,e_r)} \cap (A^+ \cdot \{(f,q_i):f \in
F\} \times A^+) \delta_A\ (i = 1, \ldots, m)$$ then we can write
$$L_{(f,e_r)} = \bigcup_{i=1,\ldots,m}{L_{(f,e_r)}^{(i)}}$$ and it
suffices to prove that, for each $i \in \{1,\ldots,m\}$, the
language $L_{(f,e_r)}^{(i)}$ is regular. To achieve that, we will
use Lemma \ref{thmgsm}. We start by showing that
$$L_{(f,e_r)}^{(i)} = K_{\bar{w}} \pi_F \zeta_{\mathcal{A}}
\delta_A \cap (A^+ \cdot \{(f,q_i): f\in F\} \times A^+ \cdot
\{(f,q_i): f\in F\})\delta_A$$ where $\bar{w}$ is the word in $K$
that represents $\ls{q_i}{f} \in S^{|T|}$. Let
$$(f_{\alpha_1},e_i) \ldots
(f_{\alpha_{n-1}},e_i)(f_{\alpha_n},q_i), (f_{\beta_1},e_j) \ldots
(f_{\beta_{s-1}},e_j)(f_{\beta_s},q_j) \in L.$$ Then
$$
\begin{array}{ll}
& ((f_{\alpha_1},e_i) \ldots
(f_{\alpha_{n-1}},e_i)(f_{\alpha_n},q_i), (f_{\beta_1},e_j) \ldots
(f_{\beta_{s-1}},e_j)(f_{\beta_s},q_j))\delta_A \in
L_{(f,e_r)}^{(i)} \\
\Leftrightarrow\; & f_{\alpha_1}\ldots f_{\alpha_{n}}
\ls{q_i}{f} = f_{\beta_1} \ldots f_{\beta_s}\ \&\
e_i q_i e_r = e_j q_j \\
\Leftrightarrow\; & f_{\alpha_1}\ldots f_{\alpha_{n}}
\ls{q_i}{f} = f_{\beta_1} \ldots f_{\beta_s}\ \&\
e_i q_i = e_j q_j \\
\Leftrightarrow\; & f_{\alpha_1}\ldots f_{\alpha_{n}}
\ls{q_i}{f} = f_{\beta_1} \ldots f_{\beta_s}\ \&\
t_i = t_j\\
\Leftrightarrow\; & (f_{\alpha_1}\ldots
f_{\alpha_{n}},f_{\beta_1} \ldots f_{\beta_s})\delta_F \in
K_{\bar{w}} \ \&\ i = j \\
\Leftrightarrow\; & ((f_{\alpha_1},e_i) \ldots
(f_{\alpha_{n-1}},e_i)(f_{\alpha_n},q_i), (f_{\beta_1},e_j) \ldots
(f_{\beta_{s-1}},e_j)(f_{\beta_s},q_j))\delta_A \in \\
& K_{\bar{w}} \pi_F \zeta_{\mathcal{A}} \delta_A \cap (A^+ \cdot
\{(f,q_i): f\in F\} \times A^+ \cdot \{(f,q_i): f\in F\})\delta_A
\end{array}
$$
We conclude, by Lemma \ref{thmgsm}, that $L_{(f,e_r)}^{(i)}$ is
a regular language. For a generator $(f,q_r) \in A$ will we prove
that $L_{(f,q_r)}$ is regular in a similar way. We can write
$$L_{(f,q_r)} = \bigcup_{i=1,\ldots,m}{L_{(f,q_r)}^{(i)}}$$ where
$$L_{(f,q_r)}^{(i)} = L_{(f,q_r)} \cap (A^+ \cdot \{(f,q_i):f \in
F\} \times A^+) \delta_A\ (i = 1,\ldots ,m).$$ We let $i \in
\{1,\ldots, m\}$ arbitrary and we will prove that
$L_{(f,q_r)}^{(i)}$ is a regular language. Let $j$ the unique
element in $\{1,\ldots, m\}$ such that $e_i q_i q_r = e_j q_j$ and
let $\bar{w}$ be the word in $K$ that represents $\ls{q_i}{f} \in
S^{|T|}$. Let
$$(f_{\alpha_1},e_i) \ldots
(f_{\alpha_{n-1}},e_i)(f_{\alpha_n},q_i), (f_{\beta_1},e_k) \ldots
(f_{\beta_{s-1}},e_k)(f_{\beta_s},q_k) \in L.$$ Then
$$
\begin{array}{ll}
& ((f_{\alpha_1},e_i) \ldots
(f_{\alpha_{n-1}},e_i)(f_{\alpha_n},q_i), (f_{\beta_1},e_k) \ldots
(f_{\beta_{s-1}},e_k)(f_{\beta_s},q_k))\delta_A \in
L_{(f,q_r)}^{(i)} \\
\Leftrightarrow\; & f_{\alpha_1}\ldots f_{\alpha_{n}}
\ls{q_i}{f} = f_{\beta_1} \ldots f_{\beta_s}\ \&\
e_i q_i q_r = e_k q_k \\
\Leftrightarrow\; & (f_{\alpha_1}\ldots
f_{\alpha_{n}},f_{\beta_1} \ldots f_{\beta_s})\delta_F \in
K_{\bar{w}} \ \&\ e_i q_i q_r = e_k q_k\ \&\ k = j \\
\Leftrightarrow\; & ((f_{\alpha_1},e_i) \ldots
(f_{\alpha_{n-1}},e_i)(f_{\alpha_n},q_i), (f_{\beta_1},e_k) \ldots
(f_{\beta_{s-1}},e_k)(f_{\beta_s},q_k))\delta_A \in \\
& K_{\bar{w}} \pi_F \zeta_{\mathcal{A}} \delta_A \cap (A^+ \cdot
\{(f,q_i): f\in F\} \times A^+ \cdot \{(f,q_j): f\in F\})\delta_A
\end{array}
$$
We can use again Lemma \ref{thmgsm} to conclude that, for each
$i$, the language
$$L_{(f,q_r)}^{(i)}
= K_{\bar{w}} \pi_F \zeta_{\mathcal{A}} \delta_A \cap (A^+ \cdot
\{(f,q_i): f\in F\} \times A^+ \cdot \{(f,q_j): f\in
F\})\delta_A$$ is regular. \qed

In the case where the semigroups $S$ and $T$ are monoids,
necessary and sufficient conditions for the wreath product
$\wre{S}{T}$ to be finitely generated are given in \cite{FGwpM}.

\begin{prop}
\label{propwpM} Let $S$ and $T$ be monoids, and let $G$ be the
group of units of $T$. Then the wreath product $\wre{S}{T}$ is
finitely generated if and only if both $S$ and $T$ are finitely
generated, and either $S$ is trivial, or $T = VG$ for some finite
subset $V$ of $T$.
\end{prop}

By using this result, our theorem has the following consequence:

\begin{cor}
Let $S$ be an automatic monoid and $T$ be a finite monoid. Then
the wreath product $\wre{S}{T}$ is automatic.
\end{cor}
\proof We assume that $S$ is not trivial. We can apply Proposition
\ref{propwpM}, with $V = T$, and so $\wre{S}{T}$ is finitely
generated. Moreover, the three conditions in Proposition
\ref{proprob}, for the case where the diagonal $S$-act is not
finitely generated, hold trivially since $S$ and $T$ are monoids.
The proof of our theorem is based on these conditions and
therefore the wreath product $\wre{S}{T}$ is automatic. \qed

It is still an open question whether or not the wreath product
$\wre{S}{T}$ is automatic when it is finitely generated. Of
course, from the above result, it only remains to consider
the case where the diagonal $S$-act is finitely generated. In
\cite{FGwpM} and \cite{Rob} we can find some examples of wreath
products with finitely generated diagonal $S$-act which, as the
authors observe, is in some way the less common case. Another
interesting problem is that of the automaticity of the wreath
product in the case where the semigroup $T$ is also infinite. A
natural starting point here is to use Proposition \ref{propwpM}
and investigate the case where $S$ and $T$ are monoids.

\newpage
\providecommand{\bysame}{\leavevmode\hbox to3em{\hrulefill}\thinspace}
\providecommand{\MR}{\relax\ifhmode\unskip\space\fi MR }
\providecommand{\MRhref}[2]{%
  \href{http://www.ams.org/mathscinet-getitem?mr=#1}{#2}
}
\providecommand{\href}[2]{#2}


\end{document}